\documentclass[envcountsect,runningheads]{llncs}
\usepackage{url}
\usepackage{algorithm}
\usepackage{algorithmic}

\usepackage{amsmath,amssymb}

\newcommand{\ignore}[1]{}
\newcommand{\jacobi}[2]{{\genfrac{(}{)}{}{}{#1}{#2}}}
\newcommand{\ndiv}{\nmid}

\newcommand{\F}{\mathbb{F}}
\newcommand{\Z}{\mathbb{Z}}
\newcommand{\p}{\mathfrak{p}}
\renewcommand{\O}{\mathcal{O}}
\renewcommand{\L}{\mathfrak{L}}
\newcommand{\characteristic}{\operatorname{char}}

\DeclareMathOperator{\Norm}{Norm}
\DeclareMathOperator{\End}{End}
\DeclareMathOperator{\disc}{disc}
\DeclareMathOperator{\Pic}{Pic}
\DeclareMathOperator{\Cl}{Cl}

\begin{document}

\title{A Subexponential Algorithm for Evaluating Large Degree
  Isogenies}
\author{David Jao \and Vladimir Soukharev}
\institute{
Department of Combinatorics and Optimization \\
University of Waterloo, Waterloo, Ontario, N2L 3G1, Canada\\
\email{\{djao,vsoukhar\}@math.uwaterloo.ca}
}

\maketitle
\begin{abstract}
  An isogeny between elliptic curves is an algebraic morphism which is
  a group homomorphism. Many applications in cryptography require
  evaluating large degree isogenies between elliptic curves
  efficiently. For ordinary curves of the same endomorphism ring, the previous
  best known algorithm has a worst case running time which is
  exponential in the length of the input. In this paper we show this
  problem can be solved in subexponential time under reasonable
  heuristics. Our approach is based on factoring the ideal
  corresponding to the kernel of the isogeny, modulo principal ideals,
  into a product of smaller prime ideals for which the isogenies can
  be computed directly. Combined with previous work of Bostan et al.,
  our algorithm yields equations for large degree isogenies in
  quasi-optimal time given only the starting curve and the kernel.
\end{abstract}
\section{Introduction}\label{sec:intro}

A well known theorem of Tate~\cite{tate} states that two elliptic
curves defined over the same finite field $\F_q$ are isogenous (i.e.
admit an isogeny between them) if and only if they have the same
number of points over $\F_q$. Using fast point counting algorithms
such as Schoof's algorithm and others~\cite{AlgBook,schoof}, it is
very easy to check whether this condition holds, and thus whether or
not the curves are isogenous. However, constructing the actual isogeny
itself is believed to be a hard problem due to the nonconstructive
nature of Tate's theorem.  Indeed, given an ordinary curve $E/\F_q$
and an ideal of norm $n$ in the endomorphism ring, the fastest
previously known algorithm for constructing the unique (up to
isomorphism) isogeny having this ideal as kernel has a running time of
$O(n^{3+\varepsilon})$, except in a certain very small number of
special cases~\cite{BCL,galbraith,GHS}. In this paper, we present a
new probabilistic algorithm for evaluating such isogenies, which in
the vast majority of cases runs (heuristically) in subexponential time.
Specifically, we show that for ordinary curves, one can evaluate
isogenies of degree~$n$ between curves of nearly equal endomorphism ring
over $\F_q$ in time less than $L_q(\frac{1}{2},\frac{\sqrt{3}}{2})
\log(n)$, provided $n$ has no large prime divisors in common with the
endomorphism ring discriminant.  Although this running time is not
polynomial in the input
length, our algorithm is still much faster than the (exponential)
previous best known algorithm, and in practice allows for the
evaluation of isogenies of cryptographically sized degrees, some
examples of which we present here. We emphasize that, in contrast with
the previous results of Br\"oker et al.~\cite{BCL}, our algorithm is
not limited to special curves such as pairing friendly curves with
small discriminant.

If an explicit equation for the isogeny as a rational function is
desired, our approach in combination with the algorithm of Bostan et
al.~\cite{bostan} can produce the equation in time
$O(n^{1+\varepsilon})$ given $E$ and an ideal of norm
$n$, which is quasi-optimal in the
sense that (up to log factors) it is equal to the size of the output. To
our knowledge, this method is the only known algorithm for computing
rational function expressions of large degree isogenies in quasi-optimal
time in the general case, given only the starting curve and the kernel.

Apart from playing a central role in the implementation of the point
counting algorithms mentioned above, isogenies have been used in
cryptography to transfer the discrete logarithm problem from one
elliptic curve to
another~\cite{AlgBook,galbraith,GHS,jmv,weak-fields,teske}.  In many of
these applications, our algorithm cannot be used directly, since in
cryptography one is usually given two isogenous curves, rather than
one curve together with the isogeny degree. However, earlier
results~\cite{galbraith,GHS,jmv} have shown that the problem of
computing isogenies between a given pair of curves can be reduced to the
problem of computing isogenies of prime degree starting from a given
curve.  It is therefore likely that the previous
best isogeny construction algorithms in the cryptographic setting can
be improved or extended in light of the work that we present here.

\section{Background}\label{sec:background}

Let $E$ and $E'$ be elliptic curves defined over a finite field $\F_q$
of characteristic $p$. An isogeny $\phi\colon E \to E'$ defined over
$\F_q$ is a non-constant rational map defined over $\F_q$ which is
also a group homomorphism from $E(\F_q)$ to~$E'(\F_q)$.  This
definition differs slightly from the standard definition in that it
excludes constant maps~\cite[\S III.4]{silverman}. The degree of an
isogeny is its degree as a rational map, and an isogeny of degree
$\ell$ is called an $\ell$-isogeny.  Every isogeny of degree greater
than 1 can be factored into a composition of isogenies of prime degree
defined over $\bar{\F}_q$~\cite{isogeny-cycles-morain}.

For any elliptic curve $E\colon y^2 + a_1 xy + a_3 y = x^3 + a_2 x^2 +
a_4 x + a_6$ defined over $\F_q$, the Frobenius endomorphism is the
isogeny $\pi_q \colon E \to E$ of degree $q$ given by the equation
$\pi_q(x,y) = (x^q,y^q)$. The characteristic polynomial of $\pi_q$ is
$X^2 - t X + q$ where $t = q+1-\#E(\F_q)$ is the trace of $E$.

An endomorphism of $E$ is an isogeny $E \to E$
defined over the algebraic closure $\bar{\F}_q$ of $\F_q$. The
set of endomorphisms of $E$ together
with the zero map forms a ring under the operations of pointwise
addition and composition; this ring is called the endomorphism ring of
$E$ and denoted $\End(E)$. The ring $\End(E)$ is isomorphic either to
an order in a quaternion algebra or to an order in an imaginary
quadratic field~\cite[V.3.1]{silverman}; in the first case we say $E$
is supersingular and in the second case we say $E$ is ordinary.

Two elliptic curves $E$ and $E'$ defined over $\F_q$ are said to be
isogenous over $\F_q$ if there exists an isogeny $\phi\colon E \to
E'$ defined over $\F_q$.  A theorem of Tate states that two curves
$E$ and $E'$ are isogenous over $\F_q$ if and only if $\#E(\F_q)
= \#E'(\F_q)$~\cite[{\S}3]{tate}.  Since every isogeny has a dual
isogeny~\cite[III.6.1]{silverman}, the property of being isogenous
over $\F_q$ is an equivalence relation on the finite set of
$\bar{\F}_q$-isomorphism classes of elliptic curves defined over
$\F_q$.  Moreover, isomorphisms between elliptic curves can be
classified completely and computed efficiently in all
cases~\cite{galbraith}. Accordingly, we define
an isogeny class to be an equivalence class of elliptic curves, taken
up to $\bar{\F}_q$-isomorphism, under this equivalence relation.

Curves in the same isogeny class are either all supersingular or all
ordinary. The vast majority of curves are ordinary, and indeed the
number of isomorphism classes of supersingular curves is finite for
each characteristic. Also, ordinary curves form the majority of the
curves of interest in applications such as cryptography. Hence, we
assume for the remainder of this paper that we are in the
\emph{\bf{ordinary case}}.

Let $K$ denote the imaginary quadratic field containing $\End(E)$,
with maximal order $\O_K$. For any order $\O \subseteq \O_K$, the
conductor of $\O$ is defined to be the integer $[\O_K:\O]$. The field
$K$ is called the CM field of $E$.  We write $c_E$ for the conductor
of $\End(E)$ and $c_\pi$ for the conductor of $\Z[\pi_q]$.  It
follows from~\cite[\S 7]{cox} that $\End(E) = \Z + c_E \O_K$ and
$\Delta = c_E^2 \Delta_K,$ where $\Delta$ (respectively, $\Delta_K$) is the
discriminant of the imaginary quadratic order $\End(E)$ (respectively,
$\O_K$).  Furthermore, the characteristic polynomial has discriminant
$\Delta_\pi = t^2 - 4q = \disc(\Z[\pi_q]) = c_\pi^2 \Delta_K$,
with $c_\pi = c_E \cdot [\End(E):\Z[\pi_q]]$.

Following~\cite{volcano} and~\cite{galbraith}, we say that an isogeny
$\phi\colon E \to E'$ of prime degree $\ell$ defined over $\F_q$ is
``down'' if $[\End(E):\End(E')] = \ell$, ``up'' if
$[\End(E'):\End(E)] = \ell$, and ``horizontal'' if $\End(E) =
\End(E)$.  Two curves in an isogeny class are said to ``have the same
level'' if their endomorphism rings are equal. Within each isogeny
class, the property of having the same level is an equivalence
relation. A horizontal isogeny always goes between two curves of the
same level; likewise, an up isogeny enlarges the endomorphism ring and
a down isogeny reduces it.  Since there are fewer elliptic curves at
higher levels than at lower levels, the collection of elliptic curves
in an isogeny class visually resembles a ``pyramid'' or a
``volcano''~\cite{volcano}, with up isogenies ascending the structure
and down isogenies descending. If we restrict to the graph of
$\ell$-isogenies for a single $\ell$, then in general the
$\ell$-isogeny graph is disconnected, having one $\ell$-volcano for
each intermediate order $\Z[\pi_q] \subset \O \subset \O_K$ such that
$\O$ is maximal at $\ell$ (meaning $\ell \ndiv [\O_K:\O]$). The ``top
level'' of the class consists of curves $E$ with $\End(E) = \O_K$, and
the ``bottom level'' consists of curves with $\End(E) = \Z[\pi_q]$.

We say that $\ell$ is an \emph{Elkies prime}~\cite[p. 119]{bss} if
$\ell \ndiv c_E$ and $\jacobi{\Delta}{\ell} \neq -1$, or equivalently
if and only if $E$ admits a horizontal isogeny of degree $\ell$. The
number of $\ell$-isogenies of each type can easily be determined
explicitly~\cite{volcano,galbraith,kohel}. In particular, for all but
the finitely many primes $\ell$ dividing $[\O_K:\Z[\pi_q]]$, we have
that every rational $\ell$-isogeny admitted by $E$ is horizontal.

\section{The Br\"oker-Charles-Lauter algorithm}\label{sec:bcl}

Our algorithm is an extension of the algorithm developed by Br\"oker,
Charles, and Lauter~\cite{BCL} to evaluate large degree isogenies
over ordinary elliptic curves with endomorphism rings of small class
number, such as pairing-friendly curves~\cite{freeman}. In this
section we provide a summary of their results.

The following notation corresponds to that of~\cite{BCL}. Let $E/\F_q$
be an ordinary elliptic curve with endomorphism ring $\End(E)$
isomorphic to an imaginary quadratic order $\O_\Delta$ of discriminant
$\Delta < 0$. Identify $\End(E)$ with $\O_\Delta$ via the unique
isomorphism $\iota$ such that $\iota^*(x) \omega = x \omega$ for all
invariant differentials $\omega$ and all $x \in \O_\Delta$. Then every
horizontal separable isogeny on $E$ of prime degree $\ell$ corresponds
(up to isomorphism) to a unique prime ideal $\L \subset \O_\Delta$ of
norm $\ell$ for some Elkies prime $\ell$. We denote the kernel of this
isogeny by $E[\L]$. Any two distinct isomorphic horizontal isogenies
(i.e., pairs of isogenies where one is equal to the composition of the
other with an isomorphism) induce different maps on the space of
differentials of $E$, and a separable isogeny is uniquely determined
by the combination of its kernel and the induced map on the space of
differentials. A \emph{normalized} isogeny is an isogeny $\phi\colon E
\to E'$ for which
$
\phi^*(\omega_{E'}) = \omega_E
$
where $\omega_E$ denotes the invariant differential of $E$.
Algorithm~\ref{BCLAlg} (identical to Algorithm~4.1 in~\cite{BCL})
evaluates, up to automorphisms of $E$, the unique normalized
horizontal isogeny of degree $\ell$ corresponding to a given kernel
ideal $\L \subset \O_\Delta$.

The following theorem, taken verbatim from~\cite{BCL}, shows that the
running time of Algorithm~\ref{BCLAlg} is polynomial in the quantities
$\log(\ell)$, $\log(q)$, $n$, and $|\Delta|$.

\begin{theorem}\label{BCLThm}
  Let $E/\F_q$ be an ordinary elliptic curve with Frobenius $\pi_q$,
  given by a Weierstrass equation, and let $P \in E(\F_{q^n})$ be a
  point on $E$. Let $\Delta = \disc(\End(E))$ be given. Assume that
  $[\End(E) : \mathbb{Z}[\pi_q]]$ and $\#E(\F_{q^n})$ are coprime, and
  let $\L = (\ell, c+d\pi_q)$ be an $\End(E)$-ideal of prime norm
  $\ell \neq \characteristic(\F_q)$ not dividing the index $[\End(E) :
  \Z[\pi_q]]$. Algorithm~\ref{BCLAlg} computes the unique elliptic
  curve $E'$ such that there exists a normalized isogeny $\phi \colon E
  \to E'$ with kernel $E[\L]$. Furthermore, it computes the
  $x$-coordinate of $\phi(P)$ if $\End(E)$ does not equal
  $\mathbb{Z}[i]$ or $\mathbb{Z}[\zeta_3]$ and the square,
  respectively cube, of the $x$-coordinate of $\phi(P)$ otherwise.
  The running time of the algorithm is polynomial in $\log(\ell)$,
  $\log(q)$, $n$ and $|\Delta|$.
\end{theorem}

\begin{algorithm}
\caption{The Br\"oker-Charles-Lauter algorithm}
\label{BCLAlg}
\begin{algorithmic}[1]
\REQUIRE A discriminant $\Delta$, an elliptic curve $E/\F_q$ with
$\End(E) = \O_{\Delta}$ and a point $P \in E(\F_{q^n})$ such
that $[\End(E) : \Z[\pi_q]]$ and $\#E(\F_{q^n})$ are coprime,
and an $\End(E)$-ideal $\L = (\ell, c+d\pi_q)$ of prime norm $\ell
\neq \characteristic(\F_q)$ not dividing the index $[\End(E) :
\mathbb{Z}[\pi_q]]$.

\ENSURE The unique elliptic curve $E'$ admitting a normalized isogeny
$\phi\colon E \to E'$ with kernel $E[\L]$, and the
$x$-coordinate of $\phi(P)$ for $\Delta \neq -3, -4$ and the square
(resp. cube) of the $x$-coordinate otherwise.

\STATE Compute the direct sum decomposition
  $\Pic(\mathcal{O}_{\Delta}) = \bigotimes\langle[I_i]\rangle$ of
  $\Pic(\mathcal{O}_{\Delta})$ into cyclic groups generated by the
  degree 1 prime ideals $I_i$ of smallest norm that are coprime to the
  product $p\cdot \#E(\F_{q^n})\cdot [\End(E) : \mathbb{Z}[\pi_q]]$.

\STATE  Using brute force\footnotemark, find $e_1, e_2, \ldots, e_k$
       such that $[\L] = [I_1^{e_1}]\cdot [I_2^{e_2}] \cdots
        [I_k^{e_k}]$.

\STATE Find $\alpha$ (using Cornacchia's algorithm) and express $\L =
I_1^{e_1}\cdot I_2^{e_2} \cdots I_k^{e_k}\cdot (\alpha)$.

\STATE Compute a sequence of isogenies $(\phi_1, \ldots, \phi_s)$ such
that the composition $\phi_c : E \rightarrow E_c$ has kernel
$E[I_1^{e_1}\cdot I_2^{e_2} \cdots I_k^{e_k}]$ using the method
of~\cite[{\S} 3]{BCL}.

\STATE Evaluate $\phi_c(P) \in E_c(\F_{q^n})$.

\STATE Write $\alpha = (u + v\pi_q)/(zm)$. Compute the isomorphism
$\eta \colon E_c \stackrel{\sim}{\to} E'$ with
$\eta^{*}(\omega_{E'}) = (u/zm)\omega_{E_c}$. Compute $Q =
\eta(\phi_{c}(P))$.

\STATE Compute $(zm)^{-1} \bmod \#E(\F_{q^n})$,
and compute $R = ((zm)^{-1}(u + v\pi_q))(Q)$.

\STATE Put $r = x(R)^{|\mathcal{O}_\Delta|^{*}/2}$ and return $(E',
       r)$.
\end{algorithmic}
\end{algorithm}

\section{A subexponential algorithm for evaluating horizontal
  isogenies}\label{sec:alg}

As was shown in Sections~\ref{sec:background} and~\ref{sec:bcl}, any
horizontal isogeny can be expressed as a composition of prime degree
isogenies, one for each prime factor of the kernel, and any prime
degree isogeny is a composition of a normalized isogeny and an
isomorphism. Therefore, to evaluate a horizontal isogeny given its
kernel, it suffices to treat the case of horizontal normalized prime
degree isogenies.

Our objective is to evaluate the unique horizontal normalized isogeny
on a given elliptic curve $E/\F_q$ whose kernel ideal is given as $\L = (\ell,
c+d\pi_q)$, at a given point $P \in E(\F_{q^n})$, where $\ell$ is an
Elkies prime. As in~\cite{BCL}, we must also impose the additional
restriction that $\ell \ndiv [\End(E) : \mathbb{Z}[\pi_q]]$; for
Elkies primes, an equivalent restriction is that $\ell \ndiv
[\O_K:\Z[\pi_q]]$, but we retain the original formulation for
consistency with~\cite{BCL}.

In practice, one is typically given $\ell$ instead of $\L$, but since
it is easy to calculate the list of (at most two) possible primes $\L$
lying over $\ell$ (cf.~\cite{BV}), these two interpretations are for all
practical purposes equivalent, and we switch freely between them when
convenient. When $\ell$ is small, one can use modular polynomial based
techniques~\cite[\S 3.1]{BCL}, which have running time $O(\ell^3
\log(\ell)^{4+\varepsilon})$~\cite{MP}.  However, for isogeny degrees
of cryptographic size (e.g. $2^{160}$), this approach is impractical.
The Br\"oker-Charles-Lauter algorithm sidesteps this problem, by using
an alternative factorization of $\L$.  However, the running time of
Br\"oker-Charles-Lauter is polynomial in $|\Delta|$, and therefore
even this method only works for small values of $|\Delta|$. In this
section we present a modified version of the Br\"oker-Charles-Lauter
algorithm which is suitable for large values of $|\Delta|$.

\footnotetext{Br\"oker, Charles, and Lauter
  mention that this computation can be done in ``various
  ways''~\cite[p. 107]{BCL}, but the only explicit method given in~\cite{BCL}
  is brute force. The use of brute force limits the algorithm to
  elliptic curves for which $|\Delta|$ is small, such as pairing-friendly
  curves.}

We begin by giving an overview of our approach. In order to handle
large values of $|\Delta|$, there are two main problems to overcome.
One problem is that we need a fast way to produce a factorization
\begin{equation}
\L = I_1^{e_1} I_2^{e_2} \cdots I_k^{e_k} \cdot (\alpha)
\label{eq:factorization}
\end{equation}
as in lines~2 and~3 of Algorithm~\ref{BCLAlg}. The other problem is
that the exponents $e_i$ in Equation~\eqref{eq:factorization} need to
be kept small, since the running times of lines~3 and~4 of
Algorithm~\ref{BCLAlg} are proportional to $\sum_i |e_i| \Norm(I_i)^2$.  The
first problem, that of finding a factorization of $\L$, can be solved
in subexponential time using the index calculus algorithm of Hafner
and McCurley~\cite{HM} (see also~\cite[Chap. 11]{BV}). To resolve the
second problem, we turn to an idea which was first introduced by
Galbraith et. al~\cite{GHS}, and recently further refined by Bisson
and Sutherland~\cite{bisson-2009}. The idea is that, in the process of
sieving for smooth norms, one can arbitrarily restrict the input
exponent vectors to sparse vectors $(e_1, e_2, ..., e_k)$ such that
$\sum_i |e_i| N(I_i)^2$ is kept small. This restriction is implemented
in line~6 of Algorithm~\ref{FactorizeIdeal}. As in~\cite{bisson-2009},
one then assumes heuristically that the imposition of this restriction
does not affect the eventual probability of obtaining a smooth norm in
the Hafner and McCurley algorithm. Note that, unlike the input
exponents, the exponents appearing in the factorizations of the ensuing
smooth norms (that is, the values of $y_i$ in
Algorithm~\ref{FactorizeIdeal}) are always small, since the norm in
question is derived from a reduced quadratic form.

We now describe the individual components of our algorithm in detail.

\subsection{Finding a factor base}
Let $\Cl(\O_\Delta)$ denote the ideal class group of $\O_\Delta$. 
Algorithm~\ref{FactorBase} produces a factor base consisting of 
split primes in $\O_\Delta$ of norm less than some bound $N$. The
optimal value of $N$ will be determined in Section~\ref{subsec:runtime}.

\begin{algorithm}[t]
\caption{Computing a factor base}
\label{FactorBase}
\begin{algorithmic}[1]
\REQUIRE A discriminant $\Delta$, a bound $N$.
\ENSURE The set $\mathcal{I}$ consisting of split prime ideals of norm
less than $N$, together with the corresponding set $\mathcal{F}$ of
quadratic forms.
\STATE Set $\mathcal{F} \leftarrow \emptyset$.
\STATE Set $\mathcal{I} \leftarrow \emptyset$.
\STATE Find all primes $p < N$ such that $(\frac{\Delta}{p}) = 1$. Call this
set $P$. Let $k = |P|$.
\STATE For each prime $p_{i} \in P$,  find an ideal
$\p_i$ of norm
$p_{i}$ (using Cornacchia's algorithm).
\STATE For each $i$, find a quadratic form $f_i = [(p_i, b_i,
  c_i)]$ corresponding to $\p_i$ in 
$\Cl(\O_\Delta)$, using the technique of~\cite[\S 3]{S}.
\STATE Output $\mathcal{I} = \{\p_1, \p_2, \ldots, \p_k\}$ and $\mathcal{F}
= \{f_1, f_2, \ldots, f_k\}$.
\end{algorithmic}
\end{algorithm}

\begin{algorithm}[h]
\caption{``Factoring'' a prime ideal}
\label{FactorizeIdeal}
\begin{algorithmic}[1]
\REQUIRE A discriminant $\Delta$, an elliptic curve $E/\F_q$ with
$\End(E) = \O_{\Delta}$, a smoothness bound $N$, a
prime ideal $\L$ of norm $\ell$ in $\O_\Delta$, an extension degree $n$.
\ENSURE Relation of the form $\L = (\alpha) \cdot
\prod_{i=1}^{k}{I_i^{e_i}}$, where $(\alpha)$ is a fractional ideal,
$I_i$ are as in Algorithm~\ref{BCLAlg}, and $e_i > 0$ are small and
sparse.
\STATE Run Algorithm \ref{FactorBase} on input $\Delta$ and $N$ to
obtain $\mathcal{I} = \{\p_1, \p_2, \ldots, \p_k\}$ and
$\mathcal{F} = \{f_1, f_2, \ldots, f_k\}$. Discard any primes dividing
$p \cdot \#E(\F_{q^n}) \cdot [\End(E) : \Z[\pi_q]]$.
\STATE Set $p_i \leftarrow \Norm(\p_i)$. (These values are also calculated
in Algorithm~\ref{FactorBase}.)
\STATE Obtain the reduced quadratic form $[\L]$ corresponding to the
ideal class of $\L$.
\REPEAT
\FOR {$i = 1, \ldots,k$}
\STATE Pick exponents $x_i$ in the range $[0, (N/p_i)^2]$ such that 
at most $k_0$ are nonzero, where $k_0$ is a global absolute constant
(in practice, $k_0 = 3$ suffices).
\ENDFOR
\STATE Compute the reduced quadratic form $\mathfrak{a} = (a, b, c)$
for which the ideal class $[\mathfrak{a}]$ is equivalent to $[\L]\cdot
\prod_{i=1}^{k}{f_i^{x_i}}$.
\UNTIL The integer $a$ factors completely into the primes $p_i$, and
the relation derived from $[\mathfrak{a}]
= [\L]\cdot \prod_{i=1}^{k}{f_i^{x_i}}$ contains fewer than $\sqrt{\log
  (|\Delta|/3)}/z$ nonzero exponents.
\STATE Write $a = \prod_{i=1}^{k}{p_i^{u_i}}$.
\FOR {i=$1, \ldots, k$}
\STATE Using the technique of Seysen (\cite[Theorem 3.1]{S}),
determine the signs of the exponents $y_i = \pm u_i$ for which
$\mathfrak{a} = \prod_{i=1}^{k}{f_i^{y_i}}$.
\STATE Let $e_i = y_i - x_i$. (These exponents satisfy $[\L] =
\prod_{i=1}^k{f_i^{e_i}}$.)
\IF {$e_i \geq 0$}
\STATE Set $I_i \leftarrow \bar{\p}_i$
\ELSE
\STATE Set $I_i \leftarrow \p_i$
\ENDIF
\ENDFOR
\STATE Compute the principal ideal $I = \L \cdot \prod_{i=1}^k{I_i^{|e_i|}}$.
\STATE Using Cornacchia's algorithm, find a generator $\beta \in
\O_\Delta$ of $I$.
\STATE Set $m \leftarrow \prod_{i=1}^k{p_i^{|e_{i}|}}$ and $\alpha \leftarrow
\frac{\beta}{m}$.
\STATE Output $\L = (\alpha)\cdot
\bar{I}_1^{|e_1|}\cdot \bar{I}_2^{|e_2|} \cdots \bar{I}_k^{|e_k|}$. 
\end{algorithmic}
\end{algorithm}

\subsection{``Factoring'' large prime degree ideals}

Algorithm~\ref{FactorizeIdeal}, based on the algorithm of Hafner and
McCurley, takes as input a discriminant $\Delta$, a curve $E$, a prime
ideal $\L$ of prime norm $\ell$ in $\O_\Delta$, a smoothness bound
$N$, and an extension degree $n$.  It outputs a factorization
\[
\L = I_1^{e_1} I_2^{e_2} \cdots I_k^{e_k} \cdot (\alpha)
\]
as in Equation~\ref{eq:factorization}, where the $I_i$'s are as in
Algorithm~\ref{BCLAlg}, the exponents $e_i$ are positive, sparse, and
small (i.e., polynomial in $N$), and the ideal $(\alpha)$ is a
principal fractional ideal generated by $\alpha$.

\subsection{Algorithm for evaluating prime degree
  isogenies}\label{subsec:EvaluatePDI}

The overall algorithm for evaluating prime degree isogenies is given
in Algorithm~\ref{EvaluatePDI}. This algorithm is identical to
Algorithm~\ref{BCLAlg}, except that the factorization of $\L$ is
performed using Algorithm~\ref{FactorizeIdeal}. To maintain
consistency with~\cite{BCL}, we have included the quantities $\Delta$
and $\End(E)$ as part of the input to the algorithm. However, we
remark that these quantities can be computed from $E/\F_q$ in
$
L_q(\frac{1}{2},\frac{\sqrt{3}}{2})
$
operations using the algorithm of Bisson and
Sutherland~\cite{bisson-2009}, even if they are not provided as input.

\begin{algorithm}[h]
\caption{Evaluating prime degree isogenies}
\label{EvaluatePDI}
\begin{algorithmic}[1]
\REQUIRE A discriminant $\Delta$, an elliptic curve $E/\F_q$ with
$\End(E) = \O_{\Delta}$ and a point $P \in E(\F_{q^n})$ such
that $[\End(E) : \Z[\pi_q]]$ and $\#E(\F_{q^n})$ are coprime,
and an $\End(E)$-ideal $\L = (\ell, c+d\pi_q)$ of prime norm $\ell
\neq \characteristic(\F_q)$ not dividing the index $[\End(E) :
\mathbb{Z}[\pi_q]]$.

\ENSURE The unique elliptic curve $E'$ admitting a normalized isogeny
$\phi\colon E \to E'$ with kernel $E[\L]$, and the
$x$-coordinate of $\phi(P)$ for $\Delta \neq -3, -4$ and the square
(resp. cube) of the $x$-coordinate otherwise.

\STATE Choose a smoothness bound $N$ (see Section~\ref{subsec:runtime}).

\STATE Using Algorithm~\ref{FactorizeIdeal} on input $(\Delta, E, N,
\L, n)$, obtain a factorization of the form $\L = I_1^{e_1}\cdot
I_2^{e_2} \cdots I_k^{e_k}\cdot (\alpha)$.

\STATE Compute a sequence of isogenies $(\phi_1, \ldots, \phi_s)$ such
that the composition $\phi_c : E \to E_c$ has kernel
$E[I_1^{e_1}\cdot I_2^{e_2} \cdots I_k^{e_k}]$ using the method
of~\cite[{\S} 3]{BCL}.

\STATE Evaluate $\phi_c(P) \in E_c(\F_{q^n})$.

\STATE Write $\alpha = (u + v\pi_q)/(zm)$. Compute the isomorphism
$\eta \colon E_c \stackrel{\sim}{\to} E'$ with
$\eta^{*}(\omega_{E'}) = (u/zm)\omega_{E_c}$. Compute $Q =
\eta(\phi_{c}(P))$.

\STATE Compute $(zm)^{-1} \bmod \#E(\F_{q^n})$,
and compute $R = ((zm)^{-1}(u + v\pi_q))(Q)$.

\STATE Put $r = x(R)^{|\mathcal{O}_\Delta|^{*}/2}$ and return $(E',
       r)$.
\end{algorithmic}
\end{algorithm}

\subsection{Running time analysis}\label{subsec:runtime}

In this section, we determine the theoretical running time of
Algorithm~\ref{EvaluatePDI}, as well as the optimal value of the
smoothness bound $N$ to use in line 1 of the algorithm. As is typical
for subexponential time factorization algorithms involving a factor
base, these two quantities depend on each other, and hence both are
calculated simultaneously.

As in~\cite{AlgBook}, we define\footnote{The definition of
  $L_n(\alpha,c)$ in~\cite{BV} differs from that of~\cite{AlgBook} in
  the $o(1)$ term. We account for this discrepancy in our text.}
$L_n(\alpha,c)$ by
\[
L_n(\alpha,c) = O(\exp((c + o(1))(\log(n))^\alpha
(\log(\log(n)))^{1-\alpha})).
\]
The quantity $L_n(\alpha,c)$ interpolates between polynomial and
exponential size as $\alpha$ ranges from $0$ to $1$. We set $N =
L_{|\Delta|}(\frac{1}{2}, z)$ for an unspecified value of $z$, and in
the following paragraphs we determine the optimal value of $z$ which
minimizes the running time of Algorithm~\ref{EvaluatePDI}. (The fact
that $\alpha = \frac{1}{2}$ is optimal is clear from the below
analysis, as well as from prior experience with integer factorization
algorithms.) For convenience, we will abbreviate $L_{|\Delta|}(\alpha,
c)$ to $L(\alpha, c)$ throughout.

Line 2 of Algorithm~\ref{EvaluatePDI} involves running
Algorithm~\ref{FactorizeIdeal}, which in turn calls
Algorithm~\ref{FactorBase}. As it turns out,
Algorithm~\ref{FactorBase} is almost the same as Algorithm 11.1
from~\cite{BV}, which requires $L(\frac{1}{2}, z)$ time, as
shown in~\cite{BV}. The only difference is that we add an additional
step where we obtain the quadratic form corresponding to each prime
ideal in the factor base. This extra step requires
$O(\log(\Norm(I))^{1+\varepsilon})$
time for a prime ideal $I$, using Cornacchia's
Algorithm~\cite{hmw}.  Thus, the overall running time for
Algorithm~\ref{FactorBase} is bounded above by
\begin{equation*}\label{eq:runtime1}
\textstyle
L(\frac{1}{2},z)
 \cdot \log(L(\frac{1}{2}, z))^{1+\varepsilon} =
L(\frac{1}{2}, z).
\end{equation*}

Line 2 of Algorithm~\ref{FactorizeIdeal} takes $\log(\ell)$ time using
standard algorithms~\cite{cox}. The loop in lines
4--9 of Algorithm~\ref{FactorizeIdeal} is very similar to the
\textsc{FindRelation} algorithm in~\cite{bisson-2009}, except that we
only use one discriminant, and we omit the requirement that $\#R/D_1 >
\#R/D_2$ (which in any case is meaningless when there is only one
discriminant).  Needless to say, this change can only speed up the
algorithm.  Taking $\mu = \sqrt{2} z$ in~\cite[Prop.  6]{bisson-2009},
we find that the (heuristic) expected running time of the loop in
lines 4--9 of Algorithm~\ref{FactorizeIdeal} is $L(\frac{1}{2},
\frac{1}{4z})$.

The next step in Algorithm~\ref{FactorizeIdeal} having nontrivial
running time is the computation of the ideal product in line 20.  To
exponentiate an element of an arbitrary semigroup to a power $e$
requires $O(\log e)$ semigroup multiplication operations~\cite[\S
1.2]{cohen}.  To multiply two ideals $I$ and $J$ in an imaginary
quadratic order (via composition of quadratic forms) requires
$O(\max(\log(\Norm(I)), \log(\Norm(J)))^{1+\varepsilon})$ bit
operations using fast multiplication~\cite[\S 6]{schonhage}. Each of
the expressions $|I_i|^{|e_i|}$ therefore requires $O(\log |e_i|)$ ideal
multiplication operations to compute, with each individual multiplication
requiring
\begin{align*}
  O((|e_i| \log(\Norm(I_i)))^{1+\varepsilon}) =
  O\left(\left(\left(\frac{N}{p_i}\right)^2
      \log(p_i)\right)^{1+\varepsilon}\right) = O(N^{2+\varepsilon})
\end{align*}
bit operations, for a total running time of $(\log e_i)
O(N^{2+\varepsilon}) = L(\frac{1}{2},2z)$ for each~$i$. This
calculation must be performed once for each nonzero exponent $e_i$. By
line 9, the number of nonzero exponents appearing in the relation is
at most $\sqrt{\log(|\Delta|/3)}/z$, so the amount of time required to
compute all of the $|I_i|^{|e_i|}$ for all $i$ is
$(\sqrt{\log(|\Delta|/3)}/z) L(\frac{1}{2},2z) = L(\frac{1}{2},2z)$.
Afterward, the values $|I_i|^{|e_i|}$ must all be multiplied
together, a calculation which entails at most
$\sqrt{\log(|\Delta|/3)}/z$ ideal multiplications where the log-norms
of the input multiplicands are bounded above by
\[
\log \Norm(I_i^{|e_i|}) = |e_i| \log \Norm(I_i) \leq
\left(\frac{N}{p_i}\right)^2 \log p_i \leq N^2 =
\textstyle{L(\frac{1}{2}, 2z)},
\]
and thus each of the (at most) $\sqrt{\log(|\Delta|/3)}/z$
multiplications in the ensuing product can be completed in time at
most $(\sqrt{\log(|\Delta|/3)}/z) L(\frac{1}{2},2z) = L(\frac{1}{2},
2z)$.  Finally, we must multiply this end result by $\L$, an operation
which requires $O(\max(\log \ell, L(\frac{1}{2},2z))^{1+\varepsilon})$
time. All together, the running time of step 20 is $L(\frac{1}{2},2z) +
O(\max(\log \ell, L(\frac{1}{2},2z))^{1+\varepsilon}) = \max((\log
\ell)^{1+\varepsilon}, L(\frac{1}{2},2z))$, and the norm of the
resulting ideal $I$ is bounded above by $\ell \cdot
\exp(L(\frac{1}{2},2z))$.

Obtaining the generator $\beta$ of $I$ in line 21 of
Algorithm~\ref{FactorizeIdeal} using Cornacchia's algorithm requires
\[ \textstyle
O(\log(\Norm(I))^{1+\varepsilon}) = (\log \ell +
L(\frac{1}{2},2z))^{1+\varepsilon}
\]
time. We remark that finding $\beta$ given $I$ is substantially easier
than the usual Cornacchia's algorithm, which entails finding $\beta$
given only $\Norm(I)$. The usual algorithm requires finding \emph{all}
the square roots of $\Delta$ modulo $\Norm(I)$, which is very slow
when $\Norm(I)$ has a large number of prime divisors. This
time-consuming step is unnecessary when the ideal $I$ itself is given,
since the embedding of the ideal $I$ in $\End(E)$ already provides (up
to sign) the correct square root of $\Delta$ mod $I$. A detailed
description of this portion of Cornacchia's algorithm in the context
of the full algorithm, together with running time figures specific to
each sub-step, is given by Hardy et al.~\cite{hmw}; for our purposes,
the running time of a single iteration of Step 6 in~\cite[\S 4]{hmw}
is the relevant figure.  This concludes our analysis of
Algorithm~\ref{FactorizeIdeal}.

Returning to Algorithm~\ref{EvaluatePDI}, we find that (as
in~\cite{BCL}) the computation of the individual isogenies $\phi_i$ in
line 3 of Algorithm~\ref{EvaluatePDI} is limited by the time required
to compute the modular polynomials $\Phi_n(x,y)$. Using the Chinese
remainder theorem-based method of Br\"oker et al.~\cite{broker-2010},
these polynomials can be computed mod $q$ in time
$O(n^3\log^{3+\varepsilon}(n))$, and the resulting polynomials require
$O(n^2(\log^2 n + \log q))$ space.  For each ideal $I_i$, the
corresponding modular polynomial of level $p_i$ only needs to be computed
once, but the polynomial once computed must be evaluated,
differentiated, and otherwise manipulated $e_i$ times, at a cost of
$O(p_i^{2+\varepsilon})$ field operations in $\F_q$ per manipulation,
or $O(p_i^{2+\varepsilon}) (\log q)^{1+\varepsilon}$ bit operations
using fast multiplication. The total running time of line 3 is
therefore
\begin{align*}
&O(p_i^{3+\varepsilon}) + \sum_i |e_i| p_i^{2+\varepsilon} (\log
q)^{1+\varepsilon} \leq
O(N^{3+\varepsilon}) + \sum_i \left(\left(\frac{N}{p_i}\right)^2\right)
p_i^{2+\varepsilon} (\log q)^{1+\varepsilon} \\
&\leq O(N^{3+\varepsilon}) +
\frac{\sqrt{\log(|\Delta|/3)}}{z} N^{2+\varepsilon} (\log
q)^{1+\varepsilon}
= \textstyle{L(\frac{1}{2},3z) + L(\frac{1}{2}, 2z) (\log q)^{1+\varepsilon}}.
\end{align*}

Similarly, the evaluation of $\phi_c$ in line 4 requires
\[
\sum_i |e_i| p_i^{2+\varepsilon} = \textstyle{L(\frac{1}{2}, 2z)}
\]
field operations in $\F_{q^n}$, which corresponds to $L(\frac{1}{2}, 2z)
(\log q^n)^{1+\varepsilon}$ bit operations using fast multiplication.

Combining all the above quantities, we obtain a total running time of

\begin{align*}
&\textstyle L(\frac{1}{2},z) & \text{(algorithm 2)} \\
&\textstyle \quad {} + L(\frac{1}{2},\frac{1}{4z}) & \text{(lines 4--9,
  algorithm 3)} \\
&\textstyle \quad {} + \max((\log \ell)^{1+\varepsilon}, L(\frac{1}{2},2z)) &
\text{(line 20, algorithm 3)} \\
&\textstyle \quad {} + (\log \ell + L(\frac{1}{2},2z))^{1+\varepsilon} &
\text{(line 21, algorithm 3)} \\
&\textstyle \quad {} + L(\frac{1}{2},3z) + L(\frac{1}{2}, 2z) (\log
q)^{1+\varepsilon} &
\text{(line 3, algorithm 4)} \\
&\textstyle \quad {} + L(\frac{1}{2},2z) (\log q^n)^{1+\varepsilon}&
\text{(line 4, algorithm 4)}
\end{align*}
\begin{align*}
&= \textstyle L(\frac{1}{2},\frac{1}{4z}) + (\log \ell + L(\frac{1}{2},
2z))^{1+\varepsilon} + L(\frac{1}{2}, 3z) + L(\frac{1}{2}, 2z)(\log
q^n)^{1+\varepsilon}.
\end{align*}

When $|\Delta|$ is large, we may impose the reasonable assumption that
$\log(\ell) \ll L(\frac{1}{2},z)$ and $\log(q^n) \ll
L(\frac{1}{2},z)$. In this case, the running time of
Algorithm~\ref{EvaluatePDI} is dominated by the expression $
L(\frac{1}{2}, \frac{1}{4z}) + L(\frac{1}{2}, 3z), $ which attains a
minimum at $z = \frac{1}{2 \sqrt{3}}$. Taking this value of $z$, we
find that the running time of Algorithm~\ref{EvaluatePDI} is equal to
$L_{|\Delta|}(\frac{1}{2}, \frac{\sqrt{3}}{2})$. Since the maximum
value of $|\Delta| \leq |\Delta_\pi| = 4q-t^2$ is $4q$, we can
alternatively express this running time as simply $L_q(\frac{1}{2},
\frac{\sqrt{3}}{2})$.

In the general case, $\log(\ell)$ and $\log(q^n)$ might be
non-negligible compared to $L(\frac{1}{2},z)$. This can
happen in one of two ways: either $|\Delta|$ is small, or (less
likely) $\ell$ is very large and/or $n$ is large. When this happens, we can
still bound the running time of Algorithm~\ref{EvaluatePDI} by taking
$z = \frac{1}{2 \sqrt{3}}$ in the foregoing calculation, although such
a choice may fail to be optimal. We then find that the running time of
Algorithm~\ref{EvaluatePDI} is bounded above by
\[
\textstyle
(\log(\ell) + L(\frac{1}{2}, \frac{1}{\sqrt{3}}))^{1+\varepsilon} +
L(\frac{1}{2}, \frac{\sqrt{3}}{2}) + 
L(\frac{1}{2}, \frac{1}{\sqrt{3}}) (\log q^n)^{1+\varepsilon}.
\]

We summarize our results in the following theorem.

\begin{theorem}\label{thm:PDI}
  Let $E/\F_q$ be an ordinary elliptic curve with Frobenius $\pi_q$,
  given by a Weierstrass equation, and let $P \in E(\F_{q^n})$ be a
  point on $E$. Let $\Delta = \disc(\End(E))$ be given. Assume that
  $[\End(E) : \mathbb{Z}[\pi_q]]$ and $\#E(\F_{q^n})$ are coprime, and
  let $\L = (\ell, c+d\pi_q)$ be an $\End(E)$-ideal of prime norm
  $\ell \neq \characteristic(\F_q)$ not dividing the index $[\End(E) :
  \Z[\pi_q]]$. Under the heuristics of~\cite[\S 4]{bisson-2009},
  Algorithm~\ref{EvaluatePDI} computes the unique elliptic curve $E'$
  such that there exists a normalized isogeny $\phi \colon E \to E'$
  with kernel $E[\L]$. Furthermore, it computes the $x$-coordinate of
  $\phi(P)$ if $\End(E)$ does not equal $\mathbb{Z}[i]$ or
  $\mathbb{Z}[\zeta_3]$ and the square, respectively cube, of the
  $x$-coordinate of $\phi(P)$ otherwise.  The running time of the
  algorithm is bounded above by
\[
\textstyle
(\log(\ell) + L(\frac{1}{2}, \frac{1}{\sqrt{3}}))^{1+\varepsilon} +
L(\frac{1}{2}, \frac{\sqrt{3}}{2}) + 
L(\frac{1}{2}, \frac{1}{\sqrt{3}}) (\log q^n)^{1+\varepsilon}.
\]
The running time of the algorithm is subexponential in $\log |\Delta|$, and
polynomial in $\log(\ell)$, $\log(q)$, and $n$.

\end{theorem}

\section{Examples}

\subsection{Small example}\label{subsec:smallex} Let $p = 10^{10} +
19$ and let $E/\F_p$ be the curve $y^2 = x^3 + 15x + 129$. Then
$E(\F_p)$ has cardinality $10000036491 = 3\cdot 3333345497$ and trace
$t = -36471$. To avoid any bias in the selection of the prime $\ell$,
we set $\ell$ to be the smallest Elkies prime of $E$ larger than
$p/2$, namely $\ell = 5000000029$. We will evaluate the $x$-coordinate
of $\phi(P)$, where $\phi$ is an isogeny of degree $\ell$, and $P$ is
chosen arbitrarily to be the point $(5940782169, 2162385016) \in
E(\F_p)$.  We remark that, although this example is designed to be
artificially small for illustration purposes, the evaluation of this
isogeny would already be infeasible if we were using prior techniques
based on modular functions of level $\ell$.

The discriminant $\Delta$ of $E$ is $\Delta = t^2 - 4p =
-38669866235$. Set $w = \frac{1+\sqrt{\Delta}}{2}$ and $\O =
\O_\Delta$. The quadratic form $(5000000029,-2326859861,270713841)$
represents a prime ideal $\L$ of norm $\ell$, and we
show how to
calculate the isogeny $\phi$ having kernel corresponding to $E[\L]$.
Using an implementation of Algorithm~\ref{FactorizeIdeal} in
MAGMA~\cite{magma}, we find immediately the relation $ \L =
(\frac{\beta}{m}) \cdot \p_{19} \cdot \p_{31}^{24} $ where $\beta =
588048307603210005w - 235788727470005542279904$, $m = 19 \cdot
31^{24}$, $\p_{19} = (19, 2w+7)$, and $\p_{31} = (31, 2w+5)$.  Using
this factorization, we can then evaluate $\phi \colon E \to E'$ using
the latter portion of Algorithm~\ref{EvaluatePDI}. We find that $E'$
is the curve with Weierstrass equation $y^2 = x^3 + 3565469415x +
7170659769$, and $\phi(P) = (7889337683, \pm 3662693258)$. We omit the
details of these steps, since this portion of the algorithm is
identical to the algorithm of Br\"oker, Charles and Lauter, and the
necessary steps are already extensively detailed in their
article~\cite{BCL}.

We can check our computations for consistency by performing a second
computation, starting from the curve $E': y^2 = x^3 + 3565469415x +
7170659769$, the point $P' = (7889337683, 3662693258) \in E'(\F_p)$,
and the conjugate ideal $\bar{\L}$, which is represented by the
quadratic form $(5000000029,2326859861,270713841)$.  Let
$\bar{\phi}\colon E' \to E''$ denote the unique normalized isogeny
with kernel $E'[\bar{\L}]$. Up to a normalization isomorphism
$\iota\colon E \to E''$, the isogeny $\bar{\phi}$ should equal the
dual isogeny $\hat{\phi}$ of $\phi$, and the composition
$\bar{\phi}(\phi(P))$ should yield $\iota(\ell P)$. Indeed, upon
performing the computation, we find that $E''$ has equation
\[
y^2 = x^3 + (15/\ell^4)x + (129/\ell^6),
\]
which is isomorphic to $E$ via the isomorphism $\iota\colon E \to E''$
defined by $\iota(x,y) = (x/\ell^2, y/\ell^3)$, and
\[
\bar{\phi}(\phi(P)) = (3163843645, 8210361642) =
(5551543736/\ell^2, 6305164567/\ell^3),
\]
in agreement with the value of $\ell P$, which is $(5551543736, 6305164567)$.

\subsection{Medium example}\label{subsec:medex}

Let $E$ be the
ECCp-109 curve~\cite{certicomb} from the Certicom ECC
Challenge~\cite{certicoma}, with equation $y^2 = x^3 + ax + b$ over
$\F_p$ where
\begin{align*}
p &= 564538252084441556247016902735257 \\
a &= 321094768129147601892514872825668 \\
b &= 430782315140218274262276694323197
\end{align*}
As before, to avoid any bias in the choice of $\ell$, we set $\ell$ to
be the least Elkies prime greater than $p/2$, and we define $w =
\frac{1+\sqrt{\Delta}}{2}$ where $\Delta = \disc(\End(E))$. Let $\L$
be the prime ideal of norm $\ell$ in $\End(E)$ corresponding to the
reduced quadratic form $(\ell,b,c)$ of discriminant $\Delta$, where $b
= -105137660734123120905310489472471$.  For each Elkies prime $p$, let
$\p_p$ denote the unique prime ideal corresponding to the reduced
quadratic form $(p,b,c)$ where $b \geq 0$. Our smoothness bound in
this case is $N = L(\frac{1}{2},\frac{1}{2\sqrt{3}}) \approx 200$.
Using Sutherland's \texttt{smoothrelation}
package~\cite{smoothrelation}, which implements the
\textsc{FindRelation} algorithm of~\cite{bisson-2009}, one finds in a
few seconds (using an initial seed of 0) the relation $\L =
\left(\frac{\beta}{m}\right) \mathfrak{I}$, where
\begin{align*}
\mathfrak{I} &= \bar{\p}_7^{72}
\bar{\p}_{13}^{100}
\bar{\p}_{23}^{14}
\bar{\p}_{47}^{2}
\bar{\p}_{73}^{2}
\bar{\p}_{103}
\p_{179}
\p_{191} \\
m &= 7^{72} 13^{100} 23^{14} 47^2 73^2 103^1 179^1 191^1
\end{align*}
and
\begin{align*}
\beta &=
3383947601020121267815309931891893555677440374614137047492987151
\backslash \\
& 2226041731462264847144426019711849448354422205800884837 \\
&{} -1713152334033312180094376774440754045496152167352278262491589014
\backslash \\
& 097167238827239427644476075704890979685 \cdot w
\end{align*}

We find that the codomain $E'$ of the normalized isogeny $\phi\colon E
\to E'$ of kernel $E[\L]$ has equation $y^2 = x^3 + a' x + b'$ where
\begin{align*}
a' &= 84081262962164770032033494307976 \\
b' &= 506928585427238387307510041944828
\end{align*}
and that the base point
{\small
\begin{align*}
P = (97339010987059066523156133908935, 149670372846169285760682371978898)
\end{align*}}
\!\!\! of $E$ given in the Certicom ECC challenge has image
{\small
\begin{align*}
(450689656718652268803536868496211, \pm
345608697871189839292674734567941).
\end{align*}}
\!\!\! under $\phi$. As with the first example, we checked the computation
for consistency by using the conjugate ideal.

\subsection{Large example}\label{sec:largex}

Let $E$ be the ECCp-239 curve~\cite{certicomb} from the Certicom ECC
Challenge~\cite{certicoma}. Then $E$ has equation $y^2 = x^3 + ax + b$
over $\F_p$ where
{\small
\begin{align*}
p\! &=\! 862591559561497151050143615844796924047865589835498401307522524859467869\\
a\! &=\! 820125117492400602839381236756362453725976037283079104527317913759073622\\
b\! &=\! 545482459632327583111433582031095022426858572446976004219654298705912499
\end{align*}}
\!\!\!Let $\L$ be the prime ideal whose norm is the least Elkies prime
greater than $p/2$ and whose ideal class is represented by the
quadratic form $(\ell,b,c)$ with $b \geq 0$. We have $N = L(\frac{1}{2},
\frac{1}{2 \sqrt{3}}) \approx 5000$, and one finds in a few hours
using \texttt{smoothrelation}~\cite{smoothrelation} that $\L$ is
equivalent to
\[
\mathfrak{I} = \bar{\p}_7^2 \p_{11} \p_{19} \p_{37}^2 \bar{\p}_{71}^2
\bar{\p}_{131} \p_{211} \bar{\p}_{389} \bar{\p}_{433} \bar{\p}_{467}
\bar{\p}_{859}^{18} \p_{863} \bar{\p}_{1019} \bar{\p}_{1151}
\bar{\p}_{1597} \bar{\p}_{2143}^6 \bar{\p}_{2207}^5 \bar{\p}_{3359}
\]
where each ideal $\p_p$ is represented by the reduced quadratic form $(p, b, c)$
having $b \geq 0$ (this computation can be reconstructed
with~\cite{smoothrelation} using the seed $7$). The quotient
$\L/\mathfrak{I}$ is generated by $\beta/m$ where $m =
\Norm(\mathfrak{I})$ and $\beta$ is
{\small
\begin{align*}
&\!-\!923525986803059652225406070265439117913488592374741428959120914067053307
\backslash\\
&4585317-917552768623818156695534742084359293432646189962935478129227909 w.
\end{align*}}
\!\!\!Given this relation, evaluating isogenies of degree
$\ell$ is a tedious but routine computation using Elkies-Atkin
techniques~\cite[\S 3.1]{BCL}. Although we do not complete it here, the
computation is well within the reach of present
technology; indeed, Br\"oker et al.~\cite{broker-2010} have computed
classical modular polynomials mod $p$ of level up to $20000$, well
beyond the largest prime of $3389$ appearing in our relation.

\section{Related work}

Bisson and Sutherland~\cite{bisson-2009} have developed an algorithm
to compute the endomorphism ring of an elliptic curve in
subexponential time, using relation-finding techniques which largely
overlap with ours. Although our main results were obtained independently,
we have incorporated their ideas into our
algorithm in several places, resulting in a simpler presentation as
well as a large speedup compared to the original version of our work.

Given two elliptic curves $E$ and $E'$ over $\F_q$ admitting a
normalized isogeny $\phi\colon E \to E'$ of degree $\ell$, the
equation of $\phi$ as a rational function contains $O(\ell)$
coefficients. Bostan et al.~\cite{bostan} have published an algorithm
which produces this equation, given $E$, $E'$, and $\ell$. Their
algorithm has running time $O(\ell^{1+\varepsilon})$, which is
quasi-optimal given the size of the output. Using our algorithm, it is
possible to compute $E'$ from $E$ and $\ell$ in time $\log(\ell)
L_{|\Delta|}(\frac{1}{2}, \frac{\sqrt{3}}{2})$ for large $\ell$.
Hence the combination of the two algorithms can produce the equation
of $\phi$ within a quasi-optimal running time of
$O(\ell^{1+\varepsilon})$, given only $E$ and $\ell$ (or $E$ and
$\L$), without the need to provide $E'$ in the input.

\section{Acknowledgments}

We thank the anonymous referees for numerous suggestions which led to
substantial improvements in our main result.

\bibliographystyle{plain}
\bibliography{paper6}
\end{document}